\documentclass{amsart}
\usepackage[utf8]{inputenc}
\pdfoutput=1
\usepackage{verbatim}
\usepackage[capitalise]{cleveref}
\usepackage{enumitem}
\usepackage{comment}
\setenumerate{label=(\roman*)}
\usepackage{tikz-cd}
\usepackage{multirow}
\usepackage{longtable}
\usepackage{nicematrix}
\usepackage{array}
\usepackage{graphicx}

\DeclareMathOperator{\im}{Im}

\DeclareMathOperator{\ann}{ann}

\newcommand{\C}{\mathbb{C}}

\renewcommand{\a}{\alpha}
\newcommand{\A}{\tilde{\alpha}}

\newcommand{\g}{\mathfrak{g}}

\newcommand{\n}{\mathfrak{n}}
\newcommand{\h}{\mathfrak{h}}

\newtheorem{theorem}{Theorem}
\newtheorem{coro}[theorem]{Corollary}
\newtheorem{prop}{Proposition}[section]
\crefname{prop}{Proposition}{Propositions}
\newtheorem{lemma}[prop]{Lemma}

\newtheorem*{lemma*}{Lemma}
\theoremstyle{definition}

\newtheorem{theoremalt}{Theorem}[theorem]

\newenvironment{theoremp}[1]{
  \renewcommand\thetheoremalt{#1}
  \theoremalt
}{\endtheoremalt}

\newtheorem{propalt}{Proposition}[prop]

\newenvironment{propp}[1]{
  \renewcommand\thepropalt{#1}
  \propalt
}{\endpropalt}

\title{Dynkin abelianisations of Lie algebras of exceptional type}
\author{Shreepranav Varma Enugandla}
\address{Indian Institute of Science, Bangalore, India}
\email{shreepranave@iisc.ac.in}

\begin{document}

\maketitle

\begin{abstract}
	We consider degenerations of all simple Lie algebras of exceptional type obtained by embedding into affine Lie algebras. We give a filtration to consider this as an abelianisation of the original Lie algebra. We then show that the associated graded of simple, finite-dimensional modules are isomorphic to Demazure modules.
\end{abstract}

\section{Introduction}

Let $\g$ be a complex simple, finite dimensional Lie algebra, and let $V(\lambda)$ be an irreducible highest weight module for $\g$. Using the PBW filtration on $V(\lambda)$, we can consider the degenerate module $\text{gr }V(\lambda)$. These degenerations have been studied extensively, especially since PBW degenerate flag varieties were defined by Feigin in \cite{Fei}. These are defined as the highest weight orbits of actions of degenerate groups on the degenerate modules $\text{gr }V(\lambda)$.

Cerulli-Irelli and Lanini showed in \cite{CL15} that these PBW degenerate flag varieties for types $A$ and $C$ are in fact Schubert varieties in a larger partial flag variety. In \cite{CLL}, the authors explore this connection further by showing that the degenerate modules in types $A$ and $C$ are isomorphic to Demazure modules for group schemes of same type but larger rank.

Properties of degenerate flag varieties are studied in detail in \cite{CIFFFR} and \cite{CIFFFR2} using methods of representation theory of quivers. Lanini and Strickland also study the cohomology of PBW degenerated flag varieties in \cite{LaSt}. One of the achievements of this framework is finding new monomial bases for $\text{gr }V(\lambda)$, and hence for $V(\lambda)$ in types $A$ and $C$ \cite{FFL11a, FFL11b}, and in a similar way for type B \cite{Mak}. The compatibility of these bases with Demazure modules for triangular Weyl group elements has been explored by Fourier in \cite{F} for type $A$, and this was extended to type $C$ by Balla, Fourier and Kambaso in \cite{BFK}. The existence of such a monomial basis for the exceptional Lie algebra $G_2$ was exhibited by Gornitskii in \cite{G}, motivating us to study flag degenerations for exceptional type Lie algebras.

In this paper, we discuss degenerations obtained by setting some commutator brackets of negative root vectors to zero in finite dimensional simple Lie algebras of exceptional type. We degenerate only the negative part since for the geometric interpretations, we consider degenerations of the highest weight modules. The goal of this paper is to study these abelianised Lie algebras and the corresponding abelianised irreducible highest weight modules.

The abelianisations we will construct are inspired by the Dynkin abelianisations for classical types described in \cite{paper}. The idea  is to embed the Dynkin diagram of a simple Lie algebra into the Dynkin diagram of a larger simple finite dimensional Lie algebra by "stretching" the diagram. Via this stretching procedure, we can make two originally adjacent nodes no longer adjacent, ensuring the corresponding root vectors commute. The easiest picture to keep in mind is that of stretching the Dynkin diagram of $A_2$ to $A_3$.

\begin{figure}[htp]
    \centering
    \includegraphics[width=4cm]{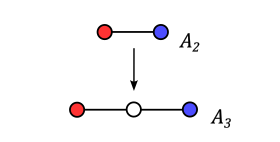}
    \label{fig:a2}
\end{figure}

In this paper, we generalize Dynkin abelianisations to the case where the larger Lie algebra is affine. This allows us to deal with the simple Lie algebras of exceptional type, which cannot be properly embedded into a Dynkin diagram of finite type. For example, the Dynkin diagram of $G_2$ can be stretched in either direction to obtain a Dynkin diagram of affine type.

\begin{figure}[htp]
    \centering
    \includegraphics[width=10cm]{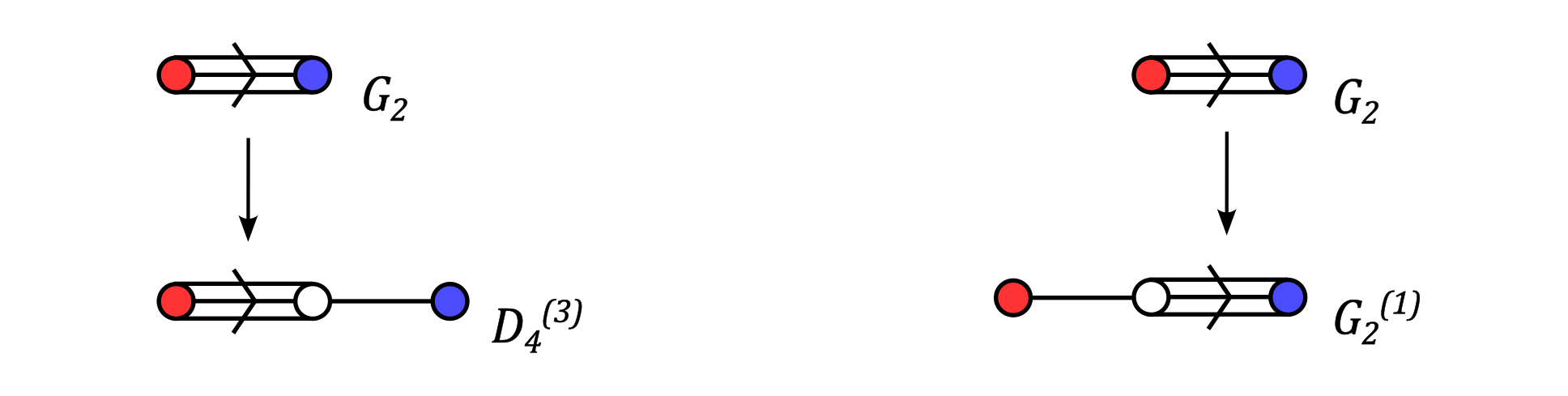}
    \label{fig:g2}
\end{figure}

For simplicity, we only deal with stretchings that make a single pair of adjacent simple roots commute. The aim is to first show that these embeddings can also be lifted to embeddings of Lie algebras, and that the corresponding abelianised modules are isomorphic to Demazure modules. Following the ideas of \cite{paper}, we first show that there is a surjection from a Demazure module to the abelianised module.


\begin{theoremp}{A}\label{thm:m}
    There exists a surjection from a Demazure module of $\hat\g$ to the abelianised module. Furthermore, if this surjection is an isomorphism for all fundamental weights, then the abeliansed module is isomorphic to the Demazure module for all weights.
\end{theoremp}

To show that this surjection is an isomorphism for all weights, it then suffices to compare dimensions for just the fundamental weights. For the cases when $\g$ is of exceptional type, there are only finitely many fundamental weights. Hence, we can check all cases by using a computer program. Once these computations are done and the surjections are verified to be isomorphisms, we obtain the main result of this paper.



\begin{propp}{B}
    For each exceptional simple Lie algebra, there exists an embedding $\iota$ into a larger (possibly affine) Lie algebra such that the abelianised modules are isomorphic to Demazure modules.
\end{propp}

The outline of this article is as follows: In \cref{sec:notation} we define the basic notations involved. In \cref{sec:main} we state the main theorem in a more precise way and show how the arguments in \cite{paper} hold even in our setting. In \cref{sec:pattern}, the embeddings are explicitly described in detail for all Lie algebras of exceptional type.\\

\textbf{Acknowledgements.} The author is supported by the Deutscher Akademischer Austauschdienst (DAAD, German Academic Exchange Service) WISE 2022 scholarship program. This project has been developed during a research visit at the Chair of Algebra and Representation Theory, RWTH Aachen University. The questions addressed in this paper were raised in a project proposal of the SFB 195.

\section{Notations}\label{sec:notation}

We follow the notation of Bourbaki \cite{bourbaki} for indexing all Dynkin diagrams.

\subsection{}
Following the notation of \cite{paper}, let $\g$ be a complex simple finite dimensional Lie algebra and fix a Cartan decomposition $\g = \n^- \oplus \h \oplus \n^+$. We consider partial degenerations of $\n^-$ constructed as follows. Fix a complex simple Lie algebra $\hat{\g}$ of finite or affine type with Cartan decomposition $\hat{\g} = \hat{\n}^- \oplus \hat{\h} \oplus \hat{\n}^+$. Let $W$ and $\hat{W}$ denote the Weyl groups corresponding to $\g$ and $\hat\g$ respectively. Finally, let $\Lambda^+$ and $\hat{\Lambda}^+$ denote the set of dominant integral weights of $\g$ and $\hat{\g}$. Let $\omega_i$ be the $i$-th fundamental weight in $\Lambda^+$, so that $\langle \omega_i, \a_j^{\vee}\rangle = \delta_{ij}$, where $\delta_{ij}$ is the Kronecker delta.\\

\subsection{}

Let $\psi:\Phi^+\to\hat\Phi^+$ be an injective function, where $\Phi^+$ and $\hat\Phi^+$ denote the sets of positive roots of $\g$ and $\hat\g$ respectively, such that $\psi$ preserves root lengths and $\im(\psi)$ is additively closed, i.e., if roots $\tilde\alpha,\tilde\beta\in \im(\psi)$ and $\tilde\a +\tilde\beta$ is a root then $\tilde\a+\tilde\beta\in \im(\psi)$. Let $\iota:\n^-\to\hat\n^-$ denote the induced map $\iota(f_\beta)=f_{\psi(\beta)}$.

$\n^{-,\iota}$ is defined to be the underlying vector space of $\n^-$ equipped with the bracket structure
\[[f_\alpha, f_\beta]_{\n^{-,\iota}} = \iota^{-1}([\iota(f_\a),\iota(f_\beta)])\]
Since $\im(\psi)$ is additively closed, it is clear that this bracket is well defined and gives a Lie algebra structure to $\n^{-,\iota}$.

We make the following assumptions on $\iota$:
\begin{itemize}
    \item $\n^{-,\iota}$ is an abelianisation of $\n^-$, i.e., for any $\a,\beta\in \Phi$, either $[f_\alpha, f_\beta]_{\n^{-,\iota}} = 0$ or $[f_\alpha, f_\beta]_{\n^{-,\iota}} = [f_\alpha, f_\beta]_{\n^-}$
    \item There exists a Weyl group element $w\in\hat{W}$ such that
    \[\{\hat\beta\in\hat\Phi^+:w(\hat\beta)\in -\hat\Phi^+\}=\im{\psi}\]
    \item There exists a monoid homomorphism $\Psi:\Lambda^+\to\hat\Lambda^+$ such that
    \[\langle \Psi(\lambda), \psi(\beta)^\vee\rangle = \langle \lambda, \beta^\vee\rangle\]
    for all $\lambda\in\Lambda^+, \beta\in\Phi^+$
\end{itemize}

For every highest weight $\g$-module $V(\lambda)$ with highest weight vector $v_\lambda$ we set $V(\lambda)^{\iota} := U(\n^{-,\iota}).v_\lambda$. We denote highest weight modules of $\hat{\g}$ by $\hat{V}(\hat{\lambda})$.

\section{The Main Theorem}\label{sec:main}

With this notation, we can restate \cref{thm:m} in a more precise way.

\begin{theorem}\label{thm:main}
    There exists a surjection $\hat{V}(\Psi(\lambda)) \supseteq U(\im\iota).v_{\Psi(\lambda)} \to V(\lambda)^\iota$ of $U(\n^{-,\iota})$-modules.
    
    Furthermore, if these surjections are isomorphisms for all fundamental weights, we have
    \[V(\lambda)^\iota \simeq U(\im\iota).v_{\Psi(\lambda)} \subseteq \hat{V}(\Psi(\lambda))\]
    for all $\lambda \in \Lambda^+$.
\end{theorem}

To prove this theorem, we need the following two lemmas. Following \cite{paper}, we let $I(\infty)\subseteq U(\n^-)$ be the left $U(\n^-)$-ideal generated by $f_{\hat\beta}$ for all $\hat\beta\in\hat\Phi^+\setminus\im\psi$.

\begin{lemma}\label{lemma:infinity}
    $U(\n^-) = U(\im\iota)\oplus I(\infty)$.
\end{lemma}

\begin{lemma}\label{lemma:ann}
    The annihilator of $v_{\Psi(\lambda)}$ in $U(\im\iota)$ is the left $U(\im\iota)$-ideal generated by all $f_\psi(\beta)^{\langle \Psi(\lambda),\psi(\beta)^\vee\rangle+1}$ for $\beta\in\im\psi$, i.e.
    \[\ann_{U(\im\iota)}v_{\Psi(\lambda)} = U(\im\iota)\langle f_{\psi(\beta)}^{\langle\Psi(\lambda),\psi(\beta)^\vee\rangle+1} \mid \beta \in \Phi^+\rangle_\C.\]
\end{lemma}

The proofs of Lemma 3.1 follows from Lemma 3.3 in \cite{paper} which goes through exactly as before, since there is no need to assume finite dimensionality of $\g$ or $\hat\g$.

Similarly, the proof of Lemma 3.2 is the exact same as the proof of Lemma 3.4 in \cite{paper}. The proof goes through even in our setting where $\hat\g$ is affine, since we only consider the finite dimensional subalgebra $U(\im\iota)$ of $U(\hat\n^-)$.\\

In \cite{paper}, the proof of Theorem 1 requires the existence of a monoid homomorphism $\Psi:\Lambda^+\to\hat\Lambda^+$ such that
    \[\langle \Psi(\lambda), \psi(\beta)^\vee\rangle = \langle \lambda, \beta^\vee\rangle\]
for all $\lambda\in\Lambda^+, \beta\in\Phi^+$, which is proved there in a previous lemma. However, the existence of such a monoid homomorphism is one of our assumptions about $\iota$, and hence the proof of Theorem 1 from \cite{paper} goes through as well.

We now want to identify this submodule of $\hat{V}(\Psi(\lambda))$ as a Demazure module.

\begin{coro}\label{coro:main}
   The map $f_\beta \mapsto e_{-w_\iota\beta}$ induces a surjection $\hat{V}_{w_\iota}(\Psi(\lambda)) \to V(\lambda)^\iota$ of $U(\n^{-,\iota})$-modules. Furthermore, if this surjection is an isomorphism for all fundamental weights, we have
   \[V(\lambda)^\iota \simeq \hat{V}_{w_\iota}(\Psi(\lambda))\]
   for all $\lambda \in \Lambda^+$.
\end{coro}

The proof of the corollary goes through exactly as in \cite{paper}.

Therefore, \cref{thm:main} and \cref{coro:main} are proved.

\section{General pattern for the embeddings}\label{sec:pattern}

We shall informally describe the general recipe used to construct the desirable embeddings.

We first define $\psi$ on the simple roots by "stretching" the tail end (an $A_k$ subgraph) of the Dynkin diagram of $\g$. For example, if the subgraph is one spanned by $\a_1,\hdots,\a_k$, then we map $\a_i$ to $\A_{i-1}$ for $i\leq k$ and $\a_i$ to $\A_i$ for $i>k$.

After defining $\psi$ on the simple roots $\a_1,\hdots,\a_n$, $\psi(\sum_{i=1}^n c_i\a_i)$ is defined to be the smallest (with respect to height) root in $\hat\Phi^+$ which has coefficient of $\psi(\a_i)$ as $c_i$ for each $i$. 

For example, consider the embedding of $G_2$ into $G_2^{(1)}$ discussed in the introduction. We wish to stretch out the root $\a_1$, so we set $\psi(\a_1)=\A_0$ and $\psi(\a_2)=\A_2$. Now we define $\psi(c_1\a_1+c_2\a_2)$ as the smallest root in $D_4^{(3)}$ such that the coefficient of $\A_0$ is $c_1$ and the coefficient of $\A_2$ is $c_2$. Hence we get 
\begin{align*}
    \psi(\a_1+\a_2) &= \A_0+\A_1+\A_2\\
    \psi(\a_1+2\a_2) &= \A_0+\A_1+2\A_2\\
    \psi(\a_1+3\a_2) &= \A_0+\A_1+3\A_2\\
    \psi(2\a_1+3\a_2) &= 2\A_0+2\A_1+3\A_2
\end{align*}

This abelianisation can also be described by giving a filtration, where the root $\beta$ is given degree $h(\beta)-c_k$, where $h(\beta)$ denotes the height of $\beta$, and $c_k$ is the coefficient of $\a_k$ in $\beta$.

For such a $\psi$, for $\a,\beta\in\Phi^+$ such that $\a+\beta\in\Phi^+$ and $\psi(\a)+\psi(\beta)\in\hat\Phi^+$, it can be seen that
\[\psi(\a+\beta)=\psi(\a)+\psi(\beta)\]
Hence,
\[\iota^{-1}([\iota(f_\a),\iota(f_\beta)]) = \iota^{-1}(kf_{\psi(\a)+\psi(\beta)})=kf_{\a+\beta}\]
for some constant $k$. Further, if we take Chevalley generators and assume that the $\a$ root string through $\beta$ is the same as the $\psi(\a)$ root string through $\psi(\beta)$, we can deduce that
\[\iota^{-1}([\iota(f_\a),\iota(f_\beta)]) = [f_\a,f_\beta]\]
showing that $\iota$ in fact induces a degeneration of $\n^-$.

It can be checked for each of the abelianisations in the following table that whenever a commutator is not set to zero, the corresponding root string structures are preserved by $\psi$.

For this $\psi$, the corresponding Weyl group element $w$ is of the form $w_0s_1\hdots s_ks_0\hdots s_{k-1}$ where $w_0$ is the longest Weyl group element for $\Phi$. The monoid homomorphism $\Psi$ is of the form 
\begin{equation*}
    \Psi(\omega_i) =
    \begin{cases}
        \tilde\omega_{i-1} & \text{if } i\leq k\\
        \tilde\omega_i & \text{if } i> k
    \end{cases}
\end{equation*}

Once we have such an embedding, \cref{thm:main} guarantees a surjection $\hat{V}_{w_\iota}(\Psi(\lambda)) \to V(\lambda)^\iota$ of $U(\n^{-,\iota})$-modules. 
For the fundamental weights $\omega_i$, we compute the dimension of the Demazure module $\hat{V}_{w_\iota}(\Psi(\omega_i))$ using a computer program (see \cite{code}), and show that this equals the dimension of $V(\omega_i)$ (which is known by the Weyl dimension formula). Since the dimension of $V(\omega_i)^\iota$ equals that of $V(\omega_i)$, it will then be clear that the surjection is in fact an isomorphism for fundamental weights. By the second part of \cref{thm:main}, this is sufficient to show that the modules are isomorphic for all $\lambda\in\Lambda^+$.

In the following table, we list out all embeddings of exceptional Lie algebras $\g$ into Dynkin diagrams with one additional node. \\

\newpage
\begin{longtable}[c]{| c || c |}
 
 \hline
 \textbf{Embedding} & 
 \begin{tabular}{ m{2cm}|m{4cm}|m{2cm}|m{2cm} }
      \textbf{Simple Roots} & \textbf{Degree} & \textbf{Weyl group element} & \textbf{Weights}\\
 \end{tabular}\\
 \hline
 \endfirsthead

 \endhead

 \hline
 \endfoot

 \endlastfoot

\hline

 $G_2$ into $D_4^{(3)}$ & 
 \begin{tabular}{ m{2cm}|m{4cm}|m{2cm}|m{2cm} } 

 $\a_1\mapsto\A_2$\newline $\a_2\mapsto\A_0$ & 
 $\a_2\mapsto 1$ \newline $\beta \mapsto c_1$ & $w_0s_1s_0$ & 
 $\omega_1\mapsto\tilde\omega_2$ \newline $\omega_2\mapsto\tilde\omega_0$\\
 
 \end{tabular}\\
 \hline
 
 $G_2$ into $G_2^{(1)}$ & 
 \begin{tabular}{ m{2cm}|m{4cm}|m{2cm}|m{2cm} } 

 $\a_1\mapsto\A_0$\newline $\a_2\mapsto\A_1$ & 
 $\a_1\mapsto 1$ \newline $\beta \mapsto c_2$ & $w_0s_1s_0$ & 
 $\omega_1\mapsto\tilde\omega_0$ \newline $\omega_2\mapsto\tilde\omega_2$\\
 
 \end{tabular}\\

\hline
 $F_4$ into $E_6^{(2)}$ & 
 \begin{tabular}{ m{2cm}|m{4cm}|m{2cm}|m{2cm} } 

 $\a_1\mapsto\A_4$\newline $\a_2\mapsto\A_3$\newline $\a_3\mapsto\A_2$\newline $\a_4\mapsto\A_0$ & $\a_4\mapsto 1$ \newline $\beta \mapsto c_1+c_2+c_3$ & 
  $w_0s_1s_0$ & 
  $\omega_1\mapsto\tilde\omega_4$ \newline $\omega_2\mapsto\tilde\omega_3$ \newline $\omega_3\mapsto\tilde\omega_2$ \newline $\omega_4\mapsto\tilde\omega_0$\\

 \hline
 
 $\a_1\mapsto\A_4$\newline $\a_2\mapsto\A_3$\newline $\a_3\mapsto\A_1$\newline $\a_4\mapsto\A_0$ & 
 $\a_3\mapsto 1$ \newline $\a_3+\a_4\mapsto 2$ \newline $\beta \mapsto c_1+c_2+c_4$ & $w_0s_1s_2s_0s_1$ & 
 $\omega_1\mapsto\tilde\omega_4$ \newline $\omega_2\mapsto\tilde\omega_3$ \newline $\omega_3\mapsto\tilde\omega_1$ \newline $\omega_4\mapsto\tilde\omega_0$\\ 
 
 \end{tabular}\\
 
 \hline
 
 $F_4$ into $F_4^{(1)}$ & 
 \begin{tabular}{ m{2cm}|m{4cm}|m{2cm}|m{2cm} }

  $\a_1\mapsto\A_0$\newline $\a_2\mapsto\A_2$\newline $\a_3\mapsto\A_3$\newline $\a_4\mapsto\A_4$ & $\a_1\mapsto 1$ \newline $\beta \mapsto c_2+c_3+c_4$ & 
  $w_0s_1s_0$ & 
  $\omega_1\mapsto\tilde\omega_0$ \newline $\omega_2\mapsto\tilde\omega_2$ \newline $\omega_3\mapsto\tilde\omega_3$ \newline $\omega_4\mapsto\tilde\omega_4$\\
 
 \hline

  $\a_1\mapsto\A_0$\newline $\a_2\mapsto\A_1$\newline $\a_3\mapsto\A_3$\newline $\a_4\mapsto\A_4$ & 
 $\a_2\mapsto 1$ \newline $\a_1+\a_2\mapsto 2$ \newline $\beta \mapsto c_1+c_3+c_4$ & $w_0s_1s_2s_0s_1$ & 
 $\omega_1\mapsto\tilde\omega_0$ \newline $\omega_2\mapsto\tilde\omega_1$ \newline $\omega_3\mapsto\tilde\omega_3$ \newline $\omega_4\mapsto\tilde\omega_4$\\
 
 \end{tabular}\\
 
 \hline
 
 $E_6$ into $E_6^{(1)}$ & 
 \begin{tabular}{ m{2cm}|m{4cm}|m{2cm}|m{2cm} } 

 $\a_1\mapsto\A_1$\newline $\a_2\mapsto\A_0$\newline $\a_3\mapsto\A_3$\newline $\a_4\mapsto\A_4$\newline $\a_5\mapsto\A_5$\newline $\a_6\mapsto\A_6$ & 
 $\a_2\mapsto 1$ \newline $\beta \mapsto h(\beta) - c_2$ & 
 $w_0s_2s_0$ & 
 $\omega_1\mapsto\tilde\omega_1$ \newline $\omega_2\mapsto\tilde\omega_0$ \newline $\omega_3\mapsto\tilde\omega_3$ \newline
 $\omega_4\mapsto\tilde\omega_4$ \newline
 $\omega_5\mapsto\tilde\omega_5$ \newline
 $\omega_6\mapsto\tilde\omega_6$\\
 
 \hline
 
 $\a_1\mapsto\A_1$\newline $\a_2\mapsto\A_0$\newline $\a_3\mapsto\A_3$\newline $\a_4\mapsto\A_2$\newline $\a_5\mapsto\A_5$\newline $\a_6\mapsto\A_6$ & 
 $\a_2\mapsto 1$ \newline $\a_2+\a_4\mapsto 2$\newline $\beta \mapsto h(\beta) - c_4$ & 
 $w_0s_2s_4s_0s_2$ & 
 $\omega_1\mapsto\tilde\omega_1$ \newline $\omega_2\mapsto\tilde\omega_0$ \newline $\omega_3\mapsto\tilde\omega_3$ \newline
 $\omega_4\mapsto\tilde\omega_2$ \newline
 $\omega_5\mapsto\tilde\omega_5$ \newline
 $\omega_6\mapsto\tilde\omega_6$\\
 
 \end{tabular}\\
 \hline
 
 $E_6$ into $E_7$ & 
 \begin{tabular}{ m{2cm}|m{4cm}|m{2cm}|m{2cm} } 

 $\a_1\mapsto\a_1$\newline $\a_2\mapsto\a_2$\newline $\a_3\mapsto\a_3$\newline $\a_4\mapsto\a_4$\newline $\a_5\mapsto\a_5$\newline $\a_6\mapsto\a_7$ & 
 $\a_6\mapsto 1$ \newline $\beta \mapsto h(\beta) - c_6$ & 
 $w_0s_6s_7$ & 
 $\omega_1\mapsto\omega_1$ \newline 
 $\omega_2\mapsto\omega_2$ \newline 
 $\omega_3\mapsto\omega_3$ \newline
 $\omega_4\mapsto\omega_4$ \newline
 $\omega_5\mapsto\omega_5$ \newline
 $\omega_6\mapsto\omega_7$\\
 
 \hline
 
 $\a_1\mapsto\a_1$\newline $\a_2\mapsto\a_2$\newline $\a_3\mapsto\a_3$\newline $\a_4\mapsto\a_4$\newline $\a_5\mapsto\a_6$\newline $\a_6\mapsto\a_7$ & 
 $\a_5\mapsto 1$ \newline $\a_5+\a_6\mapsto 2$\newline $\beta \mapsto h(\beta) - c_5$ & 
 $w_0s_6s_5s_7s_6$ & 
 $\omega_1\mapsto\omega_1$ \newline 
 $\omega_2\mapsto\omega_2$ \newline 
 $\omega_3\mapsto\omega_3$ \newline
 $\omega_4\mapsto\omega_4$ \newline
 $\omega_5\mapsto\omega_6$ \newline
 $\omega_6\mapsto\omega_7$\\
 
 \hline
 
 $\a_1\mapsto\a_1$\newline $\a_2\mapsto\a_2$\newline $\a_3\mapsto\a_3$\newline $\a_4\mapsto\a_5$\newline $\a_5\mapsto\a_6$\newline $\a_6\mapsto\a_7$ & 
 $\a_4\mapsto 1$ \newline $\a_4+\a_5\mapsto 2$\newline $\a_4+\a_5+\a_6\mapsto 3$\newline $\beta \mapsto h(\beta) - c_4$ & 
 $w_0s_6s_5s_4$-\newline -$s_7s_6s_5$ & 
 $\omega_1\mapsto\omega_1$ \newline 
 $\omega_2\mapsto\omega_2$ \newline 
 $\omega_3\mapsto\omega_3$ \newline
 $\omega_4\mapsto\omega_5$ \newline
 $\omega_5\mapsto\omega_6$ \newline
 $\omega_6\mapsto\omega_7$\\
 
 \end{tabular}\\
 
 \hline
 
 $E_7$ into $E_7^{(1)}$ & 
 \begin{tabular}{ m{2cm}|m{4cm}|m{2cm}|m{2cm} } 

 $\a_1\mapsto\A_0$\newline $\a_2\mapsto\A_2$\newline $\a_3\mapsto\A_3$\newline $\a_4\mapsto\A_4$\newline $\a_5\mapsto\A_5$\newline $\a_6\mapsto\A_6$\newline $\a_7\mapsto\A_7$ & 
 $\a_1\mapsto 1$ \newline $\beta \mapsto h(\beta) - c_1$ & 
 $w_0s_1s_0$ & 
 $\omega_1\mapsto\tilde\omega_0$ \newline $\omega_2\mapsto\tilde\omega_2$ \newline $\omega_3\mapsto\tilde\omega_3$ \newline
 $\omega_4\mapsto\tilde\omega_4$ \newline
 $\omega_5\mapsto\tilde\omega_5$ \newline
 $\omega_6\mapsto\tilde\omega_6$ \newline
 $\omega_7\mapsto\tilde\omega_7$\\
 
 \hline
 
 $\a_1\mapsto\A_0$\newline $\a_2\mapsto\A_2$\newline $\a_3\mapsto\A_1$\newline $\a_4\mapsto\A_4$\newline $\a_5\mapsto\A_5$\newline $\a_6\mapsto\A_6$\newline $\a_7\mapsto\A_7$ & 
 $\a_3\mapsto 1$ \newline $\a_1+\a_3\mapsto 2$\newline $\beta \mapsto h(\beta) - c_3$ & 
 $w_0s_1s_3s_0s_1$ & 
 $\omega_1\mapsto\tilde\omega_0$ \newline $\omega_2\mapsto\tilde\omega_2$ \newline $\omega_3\mapsto\tilde\omega_1$ \newline
 $\omega_4\mapsto\tilde\omega_4$ \newline
 $\omega_5\mapsto\tilde\omega_5$ \newline
 $\omega_6\mapsto\tilde\omega_6$ \newline
 $\omega_7\mapsto\tilde\omega_7$\\

 \hline
 
 $\a_1\mapsto\A_0$\newline $\a_2\mapsto\A_2$\newline $\a_3\mapsto\A_1$\newline $\a_4\mapsto\A_3$\newline $\a_5\mapsto\A_5$\newline $\a_6\mapsto\A_6$\newline $\a_7\mapsto\A_7$ & 
 $\a_4\mapsto 1$ \newline $\a_3+\a_4\mapsto 2$\newline $\a_1+\a_3+\a_4\mapsto 3$\newline $\beta \mapsto h(\beta) - c_4$ & 
 $w_0s_1s_3s_4$-\newline -$s_0s_1s_3$ & 
 $\omega_1\mapsto\tilde\omega_0$ \newline $\omega_2\mapsto\tilde\omega_2$ \newline $\omega_3\mapsto\tilde\omega_1$ \newline
 $\omega_4\mapsto\tilde\omega_3$ \newline
 $\omega_5\mapsto\tilde\omega_5$ \newline
 $\omega_6\mapsto\tilde\omega_6$ \newline
 $\omega_7\mapsto\tilde\omega_7$\\
 
 \end{tabular}\\
 \hline
 
 $E_7$ into $E_8$ & 
 \begin{tabular}{ m{2cm}|m{4cm}|m{2cm}|m{2cm} } 

 $\a_1\mapsto\a_1$\newline $\a_2\mapsto\a_2$\newline $\a_3\mapsto\a_3$\newline $\a_4\mapsto\a_4$\newline $\a_5\mapsto\a_5$\newline $\a_6\mapsto\a_6$\newline $\a_7\mapsto\a_8$ & 
 $\a_7\mapsto 1$ \newline $\beta \mapsto h(\beta) - c_7$ & 
 $w_0s_7s_8$ & 
 $\omega_1\mapsto\omega_1$ \newline 
 $\omega_2\mapsto\omega_2$ \newline 
 $\omega_3\mapsto\omega_3$ \newline
 $\omega_4\mapsto\omega_4$ \newline
 $\omega_5\mapsto\omega_5$ \newline
 $\omega_6\mapsto\omega_6$ \newline
 $\omega_7\mapsto\omega_8$\\
 
 \hline
 
 $\a_1\mapsto\a_1$\newline $\a_2\mapsto\a_2$\newline $\a_3\mapsto\a_3$\newline $\a_4\mapsto\a_4$\newline $\a_5\mapsto\a_5$\newline $\a_6\mapsto\a_7$\newline $\a_7\mapsto\a_8$ & 
 $\a_6\mapsto 1$ \newline $\a_6+\a_7\mapsto 2$\newline $\beta \mapsto h(\beta) - c_6$ & 
 $w_0s_7s_6s_8s_7$ & 
 $\omega_1\mapsto\omega_1$ \newline 
 $\omega_2\mapsto\omega_2$ \newline 
 $\omega_3\mapsto\omega_3$ \newline
 $\omega_4\mapsto\omega_4$ \newline
 $\omega_5\mapsto\omega_5$ \newline
 $\omega_6\mapsto\omega_7$ \newline
 $\omega_7\mapsto\omega_8$\\

 \hline
 
 $\a_1\mapsto\a_1$\newline $\a_2\mapsto\a_2$\newline $\a_3\mapsto\a_3$\newline $\a_4\mapsto\a_4$\newline $\a_5\mapsto\a_6$\newline $\a_6\mapsto\a_7$\newline $\a_7\mapsto\a_8$ & 
 $\a_5\mapsto 1$ \newline $\a_5+\a_6\mapsto 2$\newline $\a_5+\a_6+\a_7\mapsto 3$\newline $\beta \mapsto h(\beta) - c_5$ & 
 $w_0s_7s_6s_5$-\newline -$s_8s_7s_6$ & 
 $\omega_1\mapsto\omega_1$ \newline 
 $\omega_2\mapsto\omega_2$ \newline 
 $\omega_3\mapsto\omega_3$ \newline
 $\omega_4\mapsto\omega_4$ \newline
 $\omega_5\mapsto\omega_6$ \newline
 $\omega_6\mapsto\omega_7$ \newline
 $\omega_7\mapsto\omega_8$\\
 
 \hline
 
 $\a_1\mapsto\a_1$\newline $\a_2\mapsto\a_2$\newline $\a_3\mapsto\a_3$\newline $\a_4\mapsto\a_5$\newline $\a_5\mapsto\a_6$\newline $\a_6\mapsto\a_7$\newline $\a_7\mapsto\a_8$ & 
 $\a_4\mapsto 1$ \newline $\a_4+\a_5\mapsto 2$\newline $\a_4+\a_5+\a_6\mapsto 3$\newline $\a_4+\a_5+\a_6+\a_7\mapsto 4$\newline $\beta \mapsto h(\beta) - c_4$ & 
 $w_0s_7s_6s_5s_4$-\newline -$s_8s_7s_6s_5$ & 
 $\omega_1\mapsto\omega_1$ \newline 
 $\omega_2\mapsto\omega_2$ \newline 
 $\omega_3\mapsto\omega_3$ \newline
 $\omega_4\mapsto\omega_5$ \newline
 $\omega_5\mapsto\omega_6$ \newline
 $\omega_6\mapsto\omega_7$ \newline
 $\omega_7\mapsto\omega_8$\\
 
 \end{tabular}\\
 
 \hline
 
 $E_8$ into $E_8^{(1)}$ & 
 \begin{tabular}{ m{2cm}|m{4cm}|m{2cm}|m{2cm} } 

 $\a_1\mapsto\A_1$\newline $\a_2\mapsto\A_2$\newline $\a_3\mapsto\A_3$\newline $\a_4\mapsto\A_4$\newline $\a_5\mapsto\A_5$\newline $\a_6\mapsto\A_6$\newline $\a_7\mapsto\A_7$\newline $\a_8\mapsto\A_0$ & 
 $\a_8\mapsto 1$ \newline $\beta \mapsto h(\beta) - c_8$ & 
 $w_0s_8s_0$ & 
 $\omega_1\mapsto\tilde\omega_1$ \newline $\omega_2\mapsto\tilde\omega_2$ \newline $\omega_3\mapsto\tilde\omega_3$ \newline
 $\omega_4\mapsto\tilde\omega_4$ \newline
 $\omega_5\mapsto\tilde\omega_5$ \newline
 $\omega_6\mapsto\tilde\omega_6$ \newline
 $\omega_7\mapsto\tilde\omega_7$ \newline
 $\omega_8\mapsto\tilde\omega_0$\\
 
 \hline
 
 $\a_1\mapsto\A_1$\newline $\a_2\mapsto\A_2$\newline $\a_3\mapsto\A_3$\newline $\a_4\mapsto\A_4$\newline $\a_5\mapsto\A_5$\newline $\a_6\mapsto\A_6$\newline $\a_7\mapsto\A_8$\newline $\a_8\mapsto\A_0$ & 
 $\a_7\mapsto 1$ \newline $\a_7+\a_8\mapsto 2$\newline $\beta \mapsto h(\beta) - c_7$ & 
 $w_0s_8s_7s_0s_8$ & 
 $\omega_1\mapsto\tilde\omega_1$ \newline $\omega_2\mapsto\tilde\omega_2$ \newline $\omega_3\mapsto\tilde\omega_3$ \newline
 $\omega_4\mapsto\tilde\omega_4$ \newline
 $\omega_5\mapsto\tilde\omega_5$ \newline
 $\omega_6\mapsto\tilde\omega_6$ \newline
 $\omega_7\mapsto\tilde\omega_8$ \newline
 $\omega_8\mapsto\tilde\omega_0$\\

 \hline
 
 $\a_1\mapsto\A_1$\newline $\a_2\mapsto\A_2$\newline $\a_3\mapsto\A_3$\newline $\a_4\mapsto\A_4$\newline $\a_5\mapsto\A_5$\newline $\a_6\mapsto\A_7$\newline $\a_7\mapsto\A_8$\newline $\a_8\mapsto\A_0$ & 
 $\a_6\mapsto 1$ \newline $\a_6+\a_7\mapsto 2$\newline $\a_6+\a_7+\a_8\mapsto 3$\newline $\beta \mapsto h(\beta) - c_6$ & 
 $w_0s_8s_7s_6$-\newline -$s_0s_8s_7$ & 
 $\omega_1\mapsto\tilde\omega_1$ \newline $\omega_2\mapsto\tilde\omega_2$ \newline $\omega_3\mapsto\tilde\omega_3$ \newline
 $\omega_4\mapsto\tilde\omega_4$ \newline
 $\omega_5\mapsto\tilde\omega_5$ \newline
 $\omega_6\mapsto\tilde\omega_7$ \newline
 $\omega_7\mapsto\tilde\omega_8$ \newline
 $\omega_8\mapsto\tilde\omega_0$\\
 
 \hline
 
 $\a_1\mapsto\A_1$\newline $\a_2\mapsto\A_2$\newline $\a_3\mapsto\A_3$\newline $\a_4\mapsto\A_4$\newline $\a_5\mapsto\A_6$\newline $\a_6\mapsto\A_7$\newline $\a_7\mapsto\A_8$\newline $\a_8\mapsto\A_0$ & 
 $\a_5\mapsto 1$ \newline $\a_5+\a_6\mapsto 2$\newline $\a_5+\a_6+\a_7\mapsto 3$\newline $\a_5+\a_6+\a_7+\a_8\mapsto 4$\newline $\beta \mapsto h(\beta) - c_5$ & 
 $w_0s_8s_7s_6s_5$-\newline -$s_0s_8s_7s_6$ & 
 $\omega_1\mapsto\tilde\omega_1$ \newline $\omega_2\mapsto\tilde\omega_2$ \newline $\omega_3\mapsto\tilde\omega_3$ \newline
 $\omega_4\mapsto\tilde\omega_4$ \newline
 $\omega_5\mapsto\tilde\omega_6$ \newline
 $\omega_6\mapsto\tilde\omega_7$ \newline
 $\omega_7\mapsto\tilde\omega_8$ \newline
 $\omega_8\mapsto\tilde\omega_0$\\

 \hline
 
 $\a_1\mapsto\A_1$\newline $\a_2\mapsto\A_2$\newline $\a_3\mapsto\A_3$\newline $\a_4\mapsto\A_5$\newline $\a_5\mapsto\A_6$\newline $\a_6\mapsto\A_7$\newline $\a_7\mapsto\A_8$\newline $\a_8\mapsto\A_0$ & 
 $\a_4\mapsto 1$ \newline $\a_4+\a_5\mapsto 2$\newline $\a_4+\a_5+\a_6\mapsto 3$\newline $\a_4+\a_5+\a_6+\a_7\mapsto 4$\newline $\a_4+\a_5+\a_6+\a_7+\a_8\mapsto 5$\newline $\beta \mapsto h(\beta) - c_4$ & 
 $w_0s_8s_7s_6s_5s_4$-\newline -$s_0s_8s_7s_6s_5$ & 
 $\omega_1\mapsto\tilde\omega_1$ \newline $\omega_2\mapsto\tilde\omega_2$ \newline $\omega_3\mapsto\tilde\omega_3$ \newline
 $\omega_4\mapsto\tilde\omega_5$ \newline
 $\omega_5\mapsto\tilde\omega_6$ \newline
 $\omega_6\mapsto\tilde\omega_7$ \newline
 $\omega_7\mapsto\tilde\omega_8$ \newline
 $\omega_8\mapsto\tilde\omega_0$\\
 
 \end{tabular}\\
 
 \hline

\end{longtable}


\begin{prop}
Each of the embeddings in the table above gives a Dynkin abelianisation of $\g$, in the sense that 
\[[f_\a,f_\beta]_{\n^{-,\iota}}=0\ \text{or }[f_\a,f_\beta]_{\n^{-,\iota}}=[f_\a,f_\beta]_{\n^-}\]
for each $\a,\beta\in\Phi$.
\end{prop}

We know that the underlying vector space of $V(\lambda)^\iota$ is the same as that of $V(\lambda)$, and hence $\dim(V(\lambda)^\iota) = \dim(V(\lambda))$, and the latter are known by the Weyl dimension formula. To calculate the dimension of $\hat{V}_{w_\iota}(\Psi(\lambda))$, we implemented the Demazure character formula on Python (see \cite{code}). Thus we can compare the corresponding dimensions and see that they are equal, proving that the surjections given by \cref{thm:main} are isomorphisms for the fundamental weights.

This leads us to the main result of this paper.

\begin{theorem}

For each of the embeddings $\iota$ in the table above, 
\[\dim(V(\omega_i)^\iota)=\dim(\hat{V}_{w_\iota}(\Psi(\omega_i)))\]
for every fundamental weight $\omega_i$. Therefore, 
\[V(\lambda)^\iota \simeq \hat{V}_{w_\iota}(\Psi(\lambda))\]
for each $\lambda\in\Lambda^+$.
\end{theorem}


\bibliographystyle{alpha}
\bibliography{main}


\end{document}